\documentclass[11pt,a4paper]{article}
\usepackage{amsfonts,amsmath,amstext,amssymb}
\usepackage{dsfont}
\usepackage{theorem}    
\usepackage{a4wide}

\let\ddpp\mathds 
\def\R{\ddpp R}     
\def\S{\ddpp S}     
\def\L{\ddpp L}     

\newtheorem{lema}{Lemma}[section]

\newtheorem{teor}[lema]{\bf Theorem}
\newtheorem{coro}[lema]{\bf Corollary}
\theorembodyfont{\rmfamily}

\newtheorem{rema}[lema]{\bf Remark}

\renewcommand{\inf}{\mathrm{Inf}}

\title{Compact spacelike surfaces whose mean curvature function satisfies a nonlinear inequality in a 3-dimensional\\ Generalized Robertson-Walker spacetime}

\author{ Alfonso Romero${}^{a*}$ and Rafael M. Rubio${}^b$\footnote{Both authors are partially supported by the
Spanish MEC-FEDER Grant MTM2010-18099. The first author is also
partially supported by the Junta de Andaluc\'{\i}a Regional Grant P09-FQM-4496.} \\[6mm]
${}^a$ Departamento de Geometr\'{\i}a y Topolog\'{\i}a, \\
Universidad de Granada, 18071 Granada, Spain \\ E-mail\textup{:
\texttt{aromero@ugr.es}}\\[3mm]
${}^b$ Departamento de Matem\'aticas, Campus de Rabanales,\\[0.5mm] Universidad de C\'ordoba, 14071 C\'ordoba, Spain,\\[0.5mm]
E-mail\textup{:\texttt{\,rmrubio@uco.es}} }

\date{}

\begin{document}

\maketitle

\thispagestyle{empty}

\begin{abstract} Spacelike surfaces in Generalized
Robertson-Walker spacetimes whose mean curvature function
satisfies a natural nonlinear inequality are analyzed. Several
uniqueness and nonexistence results for such compact spacelike 
surfaces are proved. In the nonparametric case, new Calabi-Bernstein
type problems are solved as a consequence.
\end{abstract}

\vspace{1,5mm}


\hyphenation{ma-ni-fold}

\section{Introduction}
Let $f: I \longrightarrow \R$ be a positive smooth function on an
open interval $I=]a,b[$, $-\infty \leq a < b \leq \infty$, of the
real line $\R$ and let $(F,g)$ be a $2$-dimensional Riemannian
manifold. For each $u\in C^{\infty}(\Omega)$, $\Omega$ an open
domain in $F$, such that $u(\Omega) \subset I$ and $\mid Du \mid <
f(u)$, where $\mid Du \mid$ stands for the length of the
{gradient} $Du$ of $u$, we consider the smooth function

$$H(u)=-\,
\mathrm{div}\,\left(\frac{Du}{2f(u)\sqrt{f(u)^2-\mid
Du\mid^2}}\right)$$

\vspace{5mm}

\begin{equation}\label{primera}
\;\;\quad \qquad \qquad -\frac{f'(u)}{2\sqrt{f(u)^2 -\mid
Du\mid^2}}\left(2\,+\,\frac{\mid Du\mid^2}{f(u)^2}\right),
\end{equation}

\vspace{2mm}

\noindent where div represents the divergence operator on
$(F,g)$. When $f=1$ (constant), and $u\in C^{\infty}(\Omega)$,
$\Omega$ an open domain in $F$, such that $\mid Du \mid < 1$, then
$H(u)$ takes the more familiar form

$$H(u)=-\,
\mathrm{div}\,\left(\frac{Du}{\sqrt{1-\mid Du\mid^2}}\right),$$
and, now, $H(u)$ is the mean curvature function of the graph
$\Sigma_u=\{(u(p),p) \, : \, p\in \Omega\}$ of $u$, which is
spacelike because $\mid Du\mid < 1$, in the Lorentzian product
$I\times F$  with respect to the unit normal vector field
$$N=-\frac{1}{\sqrt{1-\mid Du\mid^2}}\;\Big(\,1\,,\,Du\,\Big).$$

In order to give an interpretation of $H(u)$ in the general case,
consider $M=I \times F$ endowed with the Lorentzian metric
\begin{equation}\label{Lorentzian_metric}
\langle\quad , \quad\rangle = -\pi^{*}_I (dt^2) + f(\pi_I)^2
\pi^{*}_{F}(g),
\end{equation}
where $\pi_I$ and $\pi_{F}$ denote the projections onto the open
interval $I$ of $\R$ and $F$, respectively; $g$ is the Riemannian
metric of $F$ and $f>0$ is a smooth function on $I$. The
Lorentzian manifold $(M, \langle\; , \;\rangle)$ is a warped
product (in the terminology of \cite[p. 204]{O-N}) with base $(I,-dt^2)$,
fiber $(F,g)$ and warping function $f$. We will refer to $(M,
\langle\, , \,\rangle)$ as a ($3$-dimensional) \emph{Generalized Robertson-Walker} (GRW)
spacetime (extending the classical notion of Robertson-Walker
spacetime to the case when $(F,g)$ is not assumed to have constant
Gauss curvature), \cite{A-R-S1}. 

\vspace{1mm}

For each $u\in C^{\infty}(\Omega)$, $u(\Omega)\subset I$, the
induced metric on $\Omega$, via the graph $\big\{(u(p),p) : p \in
\Omega \big\} \subset M$, is written on $\Omega$ as follows,

\begin{equation}\label{induced_metric}
g_u = -du^2 + f(u)^2 g,
\end{equation}
which is positive definite, if and only if $u$ satisfies $\mid Du
\mid < f(u)$. In this case,
\begin{equation}\label{normal_vector_field}
N=-\,\frac{1}{f(u)\sqrt{f(u)^2-\mid D u\,\mid^2}}\,\Big(\,
f(u)^2\, , \, Du\,\Big),
\end{equation}

\vspace{3mm}

\noindent is a unit normal vector field and the function $H(u)$,
is the mean curvature, with respect to $N$, for the spacelike
graph of $u$ in $M$.

\vspace{1mm}

Our objective here is to state several uniqueness and nonexistence
results for entire solutions to the following nonlinear elliptic
problem:
$$\qquad\qquad\qquad \qquad \hspace*{3.25cm} H(u)^2 \leq \frac{f'(u)^2}{f(u)^2},\; \hspace*{3.25cm} \qquad \qquad \qquad\mathrm{(I.1)}$$
$$\qquad\qquad\qquad \qquad \hspace*{3.55cm} \mid Du\mid < f(u),\hspace*{3.55cm} \qquad \qquad \qquad\mathrm{(I.2)}$$
where, at each value $u_0$ of $u$, the right hand side of (I.1)
is, the squared mean curvature of the spacelike slice $t=u_0$.

\vspace{1mm}

This problem, has been studied by the authors in the particular case
in which $F$ is the Euclidean plane $\R^2$  \cite{Ro-Ru3}. In the present work, 
we will deal with (I) when $F$ is a compact Riemannian surface and use a 
different approach to the one in \cite{Ro-Ru3}. In fact, in that paper a local integral
estimation of the length of the gradient of the restriction of the warping function
on a (necessarily noncompact) spacelike surface is proved in the parametric case.
Then, this is used to obtain uniqueness results both in the parametric and nonparametric 
case (see Remark \ref{refundido}). Here, we will previously study compact spacelike 
surfaces $S$ in a GRW spacetime $M$ such that its mean curvature function $H$ 
satisfies the inequality
$$
\hspace*{4.85cm} H^2 \leq \frac{f'(t)^2}{f(t)^2}, \hspace*{4.85cm}
\mathrm{(\;\widetilde I\;)}
$$
on all $S$. Let us remark that if a GRW spacetime admits a compact spacelike 
surface, then its fiber is necessarily compact, in this case, the spacetime
is called \emph{spatially closed}  (Remark \ref{nuevo}). Our approach here is global,
in fact under suitable assumptions on $M$, the inequality $\mathrm{(\;\widetilde I\;)}$ 
implies, making use of (\ref{laplacian3}), the existence of a superharmonic function on $S$
whose constancy permits to classify the compact spacelike surfaces which satisfy $\mathrm{(\;\widetilde I\;)}$
in several distinguished cases. It should be pointed out that although formula (\ref{laplacian3}) may be extended to higher dimensions, its meaning changes drastically for $n\geq 3$ (Remark \ref{nuevo4}).  

\vspace{1mm}

Later, we will derive several results in the
nonparametric case, i.e., for entire solutions $u$ to (I) on
certain compact Riemannian surfaces.

\vspace{1mm}

Several comments on (I) (and $\mathrm{(\;\widetilde I\;)}$) are now
in order: {\bf (1)} This inequality means that at any value $u_0$,
of $u$, $\mid H(u_0)\mid$ does not exceed the analogous quantity
for the spacelike slice $t=u_0.$ {\bf (2)} But it is not a
comparison assumption between extrinsic quantities of two
spacelike surfaces of $M$ (the right member is, at each $u_0$, the
squared mean curvature of $t=u_0$). {\bf (3)} It has a real sense
if $f$ is not constant (if $f$ is constant, then it is clearly
equivalent to $H=0$). {\bf (4)} In any GRW spacetime,
$\mathrm{(\;\widetilde I\;)}$ is automatically satisfied by any of
its maximal surfaces, and, under some extra assumption, by any of
its complete constant mean curvature spacelike surfaces
\cite{Ro-Ru1}. {\bf (5)} The inequality may be
physically interpreted as follows: when $F$ is compact, equation
(I.2) reads $\vert Du\vert<\lambda f(u)$, where $\lambda$ is
constant whit $0<\lambda<1$. Therefore, we get $\cosh
\theta<\frac{1}{\sqrt{1-\lambda^ 2}}$, where $\theta$ is the
hyperbolic angle between $N$ and $-\partial_t$. Along $S$ there
exist two families of instantaneous observers $\mathcal{T}_p$, $p\in
S$, where $\mathcal{T}_p=-\partial_t(p)$ (the sign minus depends of
the time orientation chosen here) and the normal observers $N_p$,
$p\in S$. The quantities $\cosh \theta (p)$ and
$v(p):=(\frac{1}{\cosh\theta(p)}) N_p^F$ are respectively the
\emph{energy} and the \emph{relative velocity} that $\mathcal{T}_p$ measures for $N_p$, and
the speed satisfies $\vert v\vert = \tanh\theta$ on $S$, \cite[\,2.1.3]{Sa-Wu}. Therefore,
we have $\vert v\vert <\lambda$ and hence $\vert v\vert$ does not approach to speed of
light in vacuum $(=1)$ on $S$. On the other hand, if the GRW spacetime is
a perfect fluid, the restriction of the total energy on any
compact spacelike surface, under the inequality assumption, is
bounded from above in terms of topological and extrinsic
quantities (see Remark \ref{energy_density}).

\vspace{2mm}

Observe that the mean curvature of the spacelike slice $t=t_0$ is
the constant $-\frac{f'(t_0)}{f(t_0)}$, so in the parametric case,
is natural to wonder,

\begin{quote} {\it Under which assumptions, the spacelike slices $t=t_0$
are characterized by the inequality $\mathrm{(\;\widetilde
I\;)}$?}
\end{quote}

On the other hand, For each $t_0\in I$, the spacelike
graph defined by the constant function $u=t_0$ satisfies
inequality (I), with $H=-\frac{f'(t_0)}{f(t_0)}$, thus, it is also
natural to wonder,

\vspace{2mm}

\begin{quote}{\it Under which assumptions, the constant
functions $u=t_0$ are the only entire solutions of the inequality
{\rm (I)}?}
\end{quote}

Our main aim in this paper is to state several answers, under
suitable assumptions on the warping function $f$ with a clear
geometric meaning. In order to do that, we will work directly on
(immersed) spacelike surfaces in a 3-dimensional GRW spacetime $M$
instead of spacelike graphs. Recall that a spacelike surface in $M$ is
locally a spacelike graph and this holds globally under some extra
topological assumptions \cite[Section 3]{A-R-S1}. Let remark that the three 
dimensional GRW spacetimes may be seen as a suitable family of toy
cosmological models where to test mathematical properties with a
potential extension to realistic four dimensional spacetimes.
The assumption ``spatially closed'' has a well-know physical meaning:
the spatial universe seems to be finite for a distinguished family of
observers, which is physically reasonable at least as a first attempt to
model spacetime. In a spatially closed GRW spacetime, the spacelike slices 
are (finite) observed spatial universes. To state when spacelike slices are
the only compact spacelike surfaces which satisfies the nonlinear inequality
is just our main goal. In fact, we prove (Theorem \ref{tc1}),

\begin{quote}
{\it Let $S$ be a compact spacelike surface of a  GRW
spacetime, whose warping function satisfies $(\log f)''\leq 0$. If
the mean curvature $H$ of $S$ satisfies $\mathrm{(\;\widetilde I\;)}$,
then $S$ is a spacelike slice.}
\end{quote}

The convexity assumption of the function $-\log f$ in this result (and in those
that follow) is automatically fulfilled if the curvature of the GRW spacetime
satisfies a natural condition with a clear physical meaning (see Subsection \ref{Ecc}).

\vspace{1mm}

Observe that no prescribed behavior for the Gauss curvature $K^F$ of $F$ is 
assumed in previous results. Under a boundedness assumption of $K^F$, we prove 
(Corollary \ref{otro}),

\begin{quote}{\it Let $M$ be a  
GRW spacetime, whose warping function satisfies $(\log f)''\leq 0$
 and the Gauss curvature of its fiber
obeys $K^ F(\pi_F)\geq c f(\pi_I)^ 2$ for some constant $c>0$. If
$M$ admits a complete spacelike surface $S$ such that the
inequality $\mathrm{(\;\widetilde I\;)},$ holds on all $S$, then $M$ must be 
spatially closed and $S$ is a spacelike slice.}
\end{quote}

As a consequence, we get (Corollary \ref{cc1}),

\begin{quote}
{\it Let $M$ be a  GRW
spacetime whose warping function satisfies $(\log f)''\leq 0$. Let $S$ a compact spacelike surface in $M$
and put $\alpha=\min_S\frac{f'(t)^2}{f(t)^2}$. If the
mean curvature function $H$ of $S$ satisfies
$$H^2\leq \alpha,$$
\noindent then $S$ must be a spacelike slice. Moreover, there is no compact
spacelike surface in $M$ whose mean curvature satisfies
$H^2<\alpha$.}
\end{quote}

Moreover, we also give a nonexistence result (Theorem \ref{BM}) and a uniqueness
result for the case os a RW spacetime of non-positive constant sectional curvature
(Theorem \ref{RW}).

\vspace{1mm}

As consequence of these parametric theorems we obtain several new
uniqueness as well as a nonexistence results of Calabi-Bernstein
type (\rm{Theorems \ref{grafo1} and
\ref{grafo3}})

\begin{quote}
{\it Let $F$ be a $2$-dimensional compact Riemannian manifold and let $f:I
\longrightarrow ]0,\infty[$ be a  
smooth function such that $(\log f)''\leq 0$. Then, the only entire solutions $u:
F\longrightarrow I$ of the inequality $\mathrm{(I)}$ are the constant
functions. Moreover, there is no entire
solution of the inequality assuming that $\mathrm{(I.1)}$ holds
strictly.}
\end{quote}

Finally, we obtain the following uniqueness and nonexistence result,

\begin{quote}{\it Let $F$ be a $2$-dimensional compact Riemannian manifold 
and let $f:I \longrightarrow ]0,\infty[$ be a smooth function
with $(\log f)''\leq 0$  and such that
$\inf \,\frac{(f')^2}{f^2}=\alpha>0$. Then, the only entire
solutions $u: F\longrightarrow I$ of the inequality
$$H(u)^2\leq\alpha ,$$
$$\vert Du\vert<f(u),$$ are the constant functions. Moreover, there
is no entire solution if the first inequality is strict.}
\end{quote}

\section{Preliminaries} Let $f$ be a positive smooth function
defined on an open interval $I$ of $\R$ and let $(F,g)$ be a
Riemannian surface. Consider $M=I\times F$ endowed with the Lorentzian 
metric $\langle \; , \; \rangle$ given in (\ref{Lorentzian_metric}), i.e.,
$M$ is a GRW spacetime as defined in previous section. The 
vector field $\partial_{t} := \partial / \partial t \in \mathfrak{X}(M)$  is timelike (and unitary) 
and hence $M$ is time orientable. Along this paper, when a time orientation
in necessary to choose on $M$, we agree to consider
$M$ endowed with the time orientation defined by $-\partial_t$.  

The spacetime $M$ has a distinguished infinitesimal conformal symmetry, namely 
the vector field $\xi: = f({\pi}_I)\,\partial_t,$ is timelike and, from the
relationship between the Levi-Civita connections of $M$ and those
of the base and the fiber \cite[Cor. 7.35]{O-N}, satisfies
\begin{equation}\label{conformal_symmetry}
\overline{\nabla}_X\xi = f'({\pi}_I)\,X,
\end{equation}
for any $X\in \mathfrak{X}(M)$, where $\overline{\nabla}$ is the
Levi-Civita connection of the metric (\ref{Lorentzian_metric}). Thus, $\xi$
is conformal with $\mathcal{L}_{\xi}\langle \,,\, \rangle
=2\,f'({\pi}_I)\,\langle \,,\, \rangle$ and its metrically
equivalent 1-form is closed.
\par
\medskip
\subsection{Energy conditions}\label{Ecc}
Recall that a Lorentzian manifold obeys the \emph{timelike convergence
condition} (TCC) if its Ricci tensor $\overline{\mathrm{Ric}}$
satisfies
\[
\overline{\mathrm{Ric}}(Z,Z)\geq 0,
\]
for all timelike tangent vector $Z$. It is normally argued that TCC on a spacetime is the
mathematical way to express that gravity, on average, attracts \cite[p. 340]{O-N}. A
weaker energy condition is the \emph{null convergence condition} (NCC)
which reads
 $$\overline{\mathrm{Ric}}(Z,Z) \geq 0,$$ for any \emph{null} tangent vector
$Z$, i.e. $Z\neq 0$ satisfying $\langle Z,Z \rangle =0$. Clearly,
a continuity argument shows that TCC implies NCC. A spacetime $M$
obeys the \emph{ubiquitous energy condition} if
$$\overline{\mathrm{Ric}}(Z,Z)>0,$$ for all timelike tangent vector $Z$.
This last energy condition is clearly stronger than TCC and roughly means
a real presence of matter at any event of the spacetime. 

\vspace{1mm}

It is easy to see that a GRW  spacetime $M$ with a $2$-dimensional fiber
$(F,g)$ obeys TCC if and only if its warping function satisfies $f''\leq
0$ and the Ricci tensor of the fiber satisfies
${\mathrm{Ric}^F}\geq (ff''-f'^2)g$, \cite{A-R-S1}. Moreover, TCC holds if and only if NCC hold and
$f'' \leq 0$. On the other hand, if a GRW spacetime $M$ obeys the
ubiquitous energy condition, then $f''<0$. Observe that if
$f''\leq 0$ (resp. $f''<0$) then $(\log f)''\leq 0$ (resp. $(\log
f)''<0$).

\vspace{1mm}

From \cite[Cor. 7.43]{O-N}, we have
\begin{equation}\label{Riccibarra}
\overline {\mathrm{Ric}}(X,Y)=\mathrm{Ric}^F(X^F,Y^F)+\left(\,
\frac{f''}{f} + \frac{(f')^2}{f^2}\, \right) \langle
X^F,Y^F\rangle  -
        \frac{2f''}{f} \langle X,\partial_t\rangle \langle Y,\partial_t\rangle ,
\end{equation} for any tangent vectors $X$, $Y$ to $M$, where $Z^F:=Z+\langle Z, \partial_t\rangle\partial_t$ 
is the projection of  the tangent vector $Z$ on the fiber.
Therefore,
$$\overline{\mathrm{Ric}}(Z,Z)=\Big\lbrace\frac{K^F(\pi_{_F})}{f^2}-(\log f)''\Big\rbrace\langle Z,\partial_t\rangle^2,$$
in the case $Z$ is null,  where $K^F$ denotes the Gauss
curvature of $F$. Thus, the GRW spacetime $M$ obeys NCC if and
only if 
\begin{equation}\label{practica}
\frac{K^F(\pi_{_F})}{f^2}-(\log f)''\geq 0.
\end{equation}

Along the main results of this paper we will assume the warping function of
the GRW spacetime $M$ satisfies $(\log f) '' \leq 0$. Note that, in
particular, this holds whenever $M$ obeys TCC. Moreover, if $K^F \geq 0$, 
the inequality $(\log f) ''\leq 0$ is implies NCC. On the other hand, if $K^F\leq 0$
holds, then NCC implies $(\log f) '' \leq 0$. In addition, when the spacetime obeys the ubiquitous energy condition the inequality $(\log f)''<0$ is satisfied.

\subsection{The restriction of the warping function on a spacelike surface}
Let $x : S \longrightarrow M$ be a (connected) spacelike surface
in $M$, i.e., $x$ is an immersion and it induces a Riemannian
metric on the $2$-dimensional manifold $S$ from the Lorentzian
metric (\ref{Lorentzian_metric}).  As usual, we agree to represent the
induced metric with the same symbol as the one used in (\ref {Lorentzian_metric}).
Then the time-orientability of $M$ allows us to consider $N
\in \mathfrak{X}^{\perp}(S)$ as the only, globally defined,
unitary timelike normal vector field on $S$ in the same
time-orientation of $-\partial_{t}$. Thus, from the wrong way
Cauchy-Schwarz inequality, (see \cite[Prop. 5.30]{O-N}, for
instance) we have $\langle N,\partial_t\rangle \geq 1$ and
$\langle N,\partial_t\rangle=1$ at a point p if and only if
$N(p)=-\partial_t(p)$. 

\vspace{1mm}

In any GRW spacetime $M$, the level surfaces of the function
$\pi_I : M \longrightarrow I$ constitute a distinguished family of spacelike
surfaces, the so-called \emph{spacelike slices}. We agree to
represent by $t=t_0$ the spacelike slice $\{t_0\}\times F$. For a given spacelike surface  
$x : S \longrightarrow M$, we have that $x(S)$ is contained in $t=t_0$ if and only if
$\pi_I \circ x=t_0$ on $S$. Note that this holds if and only if the surface $S$ is
orthogonal to $\partial_t$ or, equivalently, orthogonal to $\xi$. We will say that $S$ is a spacelike slice if $x(S)$ equals to
$t=t_0$, for some $t_0\in I$

\vspace{1mm}

Denote by $\partial_{t}^{\top}:=\partial_{t}\,+\langle N,\partial_t
\rangle N$ the tangential component of $\partial_{t}$ on $S$. From (\ref{conformal_symmetry}) using
the Gauss formula, we have
\begin{equation}\label{gradient}
\nabla t= -\partial_{t}^{\top},
\end{equation}
where $\nabla t$ is the gradient of $t:= \pi_I \circ x$. Now, using again
the Gauss formula, taking into account
$\xi^{\top}=f(t)\,\partial_{t}^{\top}$ and (\ref{gradient}), the Laplacian
of  the function $t$ on $S$ satisfies
\begin{equation}\label{Laplacian}
\Delta t = - \frac{f'(t)}{f(t)}\Big\{ 2 + \mid \nabla t \mid^2
\Big\} -2 H \,\langle N,\partial_t \rangle,
\end{equation}
where $f(t):=f \circ t$, \, $f'(t):=f' \circ t$ and 
$H:=-(1/2)\,\mathrm{trace}(A)$ is called the \emph{mean curvature} function of $S$ relative to
$N$, where $A$ is the shape operator
associated to $N$. A spacelike surface $S$ with constant mean curvature is a
critical point of the area functional under a certain volume
constraint (see \cite{B-C}, for instance). A spacelike
surface with constant $H$ is called  a spacelike surface of constant mean
curvature surface (CMC). Note that, with our choice of $N$, the
shape operator of the spacelike slice $t=t_0$ is
$A=(f'(t_0)/f(t_0))\,I$ and $H=-f'(t_0)/f(t_0)$, therefore the
spacelike slices are totally umbilical CMC spacelike surfaces.

\vspace{1mm}

A direct computation from  (\ref{gradient}) and (\ref{Laplacian})
gives
\begin{equation}\label{Laplacian2}
\Delta f(t)=  - 2\,\frac{f'(t)^2}{f(t)} + f(t)(\log f)''(t)\mid
\nabla t \mid^2-2f'(t)H\langle N,\partial_t\rangle,
\end{equation}
for any spacelike surface $S$ in the GRW spacetime $M$.

\begin{rema}\label{nuevo}
For any spacelike surface $S$ the mapping $\pi_F\circ x$ is a
local diffeomorphism, which is a covering map when $S$ is complete
and $f(t)$ is bounded on $S$ \cite[Lemma 3.1]{A-R-S1}. Therefore
if $S$ is compact, then $F$ must be also compact.
\end{rema}
\begin{rema}\label{nuevo2}
As explained in the Introduction, the graph in $M$ of  $u \in C^{\infty}(F)$, 
which satisfies $u(F) \subset I$ and $|Du|<f(u)$ is a spacelike surface in $M$. 
Each spacelike surface is locally the graph of a function, but it is not true
globally in general.
However, if $F$ is assumed to be compact and simply connected, then
every compact spacelike surface in $M$ must be diffeomorphic to $F$ and
$x(S)$ may be seen as a spacelike graph in $M$, \cite[Prop.
3.3]{A-R-S1}.
\end{rema}

From (\ref{Laplacian2}) and taking into account $-1=\mid \nabla t \mid^2 - \langle N,\partial_t\rangle^2$, 
the Laplacian of the function $\log f(t)$ on $S$ satisfies

\begin{equation}\label{laplacian3}
\Delta\log f(t)=-\Big(\frac{f'(t)}{f(t)} + H\langle N,\partial_t\rangle\Big)^2 + \Big( H^2-\frac{f'(t)^2}{f(t)^2}\Big)\langle N,\partial_t\rangle^2 + (\log f)''(t)\vert\nabla t\vert^2.
\end{equation}
Observe that, under the assumption $(\log f)''(t)\leq 0$,  the inequality $\mathrm{(\;\widetilde I\;)}$
gives that $\log f(t)$ is superharmonic on $S$. In particular,  this holds on each spacelike surface 
in $M$ whose mean curvature function satisfies the inequality $\mathrm{(\;\widetilde I\;)}$, 
whenever $M$ obeys TCC or NCC with $K^F \leq 0$.

\begin{rema}\label{nuevo4}
In view of previous assertion, equation (\ref{laplacian3}) is a key fact for the obtainment of
the results in this paper. It is then natural to wonder if (\ref{laplacian3}) may be stated for a 
spacelike hypersurface in an $(n+1)$-dimensional GRW spacetime, for any
$n\geq 2$.  In this general case, it is not difficult to get, 
\[
\Delta\log f(t)=-\Big(\frac{f'(t)}{f(t)} + H\langle N,\partial_t\rangle\Big)^2 + \Big( H^2-\frac{f'(t)^2}{f(t)^2}\Big)\langle N,\partial_t\rangle^2 + (\log f)''(t)\vert\nabla t\vert^2 
\]

\[
\hspace*{-2.9cm}- \,(n-2)\,\frac{f'(t)^2}{f(t)^2} - (n-2)\,\frac{f'(t)}{f(t)}H \langle N,\partial_t\rangle.
\]
\end{rema}

\subsection{The Gauss curvature of a spacelike surface}
Denote by $R$ and $\overline R$  the curvature tensors of  a spacelike surface $S$
and of $M$, respectively. The Gauss equation reads
\begin{equation}\label{Curvatura}
\langle R(X,Y)U,V\rangle =\langle {\overline R}(X,Y)U,V\rangle -
 \langle AY,U\rangle \langle AX,V\rangle +\langle AX,U\rangle \langle AY,V\rangle,
\end{equation}
where $X,Y,U,V \in \mathfrak{X}(S)$. Moreover, we have the Codazzi
equation which, taking into account that the normal bundle of the
spacelike surface is negative definite, is written as follows
\begin{equation}\label{Codazzi}
{\overline R}(X,Y)N = -(\nabla_X A)Y+(\nabla_Y A)X,
\end{equation}
for all $X,Y \in \mathfrak{X}(S)$.

\vspace{1mm}

From the Gauss equation (\ref{Curvatura}) we get
\begin{equation}\label{Ricci}
\mathrm{Ric}(X,Y)=\overline {\mathrm{Ric}}(X,Y) + \langle\overline
{R}(N,X)Y,N \rangle +\,2H \langle AX,Y\rangle + \langle
A^2X,Y\rangle ,
\end{equation}
where $\mathrm{Ric}$ denotes the Ricci tensor of $S$.

\vspace{1mm}

Now we take a local orthonormal frame field $E_1,E_2,E_3$ on $M$
which is adapted to $S$, i.e., on $S$, $E_1,E_2$ are tangent to
$S$ and $E_3=N$. From (\ref{Ricci}) we obtain
\begin{equation}\label{K1}
2K = \sum_{i=1}^2 \overline{\mathrm{Ric}}(E_i,E_i) + \sum_{i=1}^2 \langle\overline
{R}(N,E_i)E_i,N\rangle -\,4H^2 + {\rm trace}(A^2),
\end{equation}
where $K$ is the Gauss curvature of $S$. We can rewrite
(\ref{K1}), using  \cite[Prob. 7.13]{O-N}, as follows
\begin{equation}\label{K2}
K =\frac{f'(t)^2}{f(t)^2}
+\Big\lbrace\frac{K^F(\pi_{_F})}{f(t)^2}- (\log f)''(t)\Big\rbrace
\mid \nabla t \mid^2+\frac{K^F(\pi_{_F})}{f(t)^2} - 2H^2 +
\frac{1}{2}\,\mathrm{trace}(A^2).
\end{equation}
The Schwarz inequality for symmetric operators on a $2$-dimensional euclidean 
vector space implies $H^2\leq \frac{1}{2}\,\mathrm{trace}(A^2)$. Taking this in mind,
if we assume $M$ obeys NCC and the inequality $\mathrm{(\;\widetilde I\;)}$  holds on $S$,
then (\ref{K2}) implies $K \geq \frac{K^F(\pi_F)}{f(t)^2}$, i.e.,
at each $p\in S$, $K(p)$ is at least the Gauss curvature of the slice
$t=t(p)$ at
the point $\pi_F(x(p))$.

\vspace{1mm}

\begin{rema}\label{energy_density}
Assume the GRW spacetime $M$ is a ($3$-dimensional) \emph{perfect
fluid} with flow vector field $-\partial_t$ and \emph{energy density function}
$\rho$ \cite[Def. 12.4]{O-N}. From (\ref{Riccibarra}) we have

\begin{equation}\label{energy}
8\pi\rho=\frac{K^F}{f^2}+\frac{(f')^2}{f^2}.
\end{equation}
Now consider a spacelike surface $S$ in $M$. Using that $M$ satisfies NCC, from
(\ref{K2}) and previous formula we have
\begin{equation}\label{desigualdad}
K\geq 8\pi\rho-H^2.
\end{equation}
on each spacelike surface $S$ in $M$. If $S$ is in addition assumed to be 
compact, the \emph{total energy} on $S$ is 
$$E_{_S}:=\int_S\rho\,dS.$$ 

The inequality $\mathrm{(\;\widetilde I\;)}$ permits to deduce from (\ref{desigualdad}),
via the Gauss-Bonnet theorem the following upper bound for the total energy on $S$,

\begin{equation}
E_{_S}\,\leq \,\frac{1}{2}\,\mathcal{X}(S) + \frac{1}{8\pi}\int_S
\frac{f'(t)^2}{f(t)^2}\,dS,
\end{equation}

\noindent where $\mathcal{X}(S)$ denotes the Euler number
of $S$.
\end{rema}

\section{Uniqueness results in the parametric case}
We begin with a direct consequence of formula (\ref{laplacian3}),
\begin{teor}\label{tc1}
Let $S$ be a compact spacelike surface of a  GRW spacetime,
whose warping function satisfies $(\log f)''\leq 0$. If the mean
curvature function $H$ of $S$ satisfies $\mathrm{(\;\widetilde I\;)}$,
then $S$ is a spacelike slice.
\end{teor}
\emph{Proof.} The assumptions on $f$ and $H$ clearly give $\Delta\log
f(t)\leq 0$ making use of (\ref{laplacian3}). The compactness of $S$ implies 
that the function $\log f(t)$ is constant on $S$, and, therefore, $f(t)$ is so.
Consider now a primitive function $\mathcal{F}$ of $f$ and write $\mathcal{F}(t)$ for 
the restriction of $\mathcal{F}\circ \pi_I$ on $S$. Note that $\nabla \mathcal{F}(t)=f(t)\,\nabla t$.
Observe that the vanishing of the first term on the right of (\ref{laplacian3}) means $\Delta\mathcal{F}(t)=0.$ 
Consequently, $\mathcal{F}(t)$ is constant, and then, $S$ is a spacelike slice.\hfill{$\Box$}

\begin{rema} Note that no prescribed behavior for the Gauss curvature $K^F$ of $F$ is 
assumed in Theorem \ref{tc1}. On the other hand, the assumption $(\log
f)''\leq 0$ can not be removed. To support this assertion consider $I=\R$, $f(t)=\cosh (t)$ 
and $F=\S^2$ with its canonical metric of constant Gauss curvature $1$, then the corresponding
$3$-dimensional GRW spacetime $M$ is the De Sitter spacetime
$\S_1^3$ of constant sectional curvature $1$. It is well-known that any compact maximal
surface in the $3$-dimensional De Sitter spacetime must be totally geodesic \cite{Ramanathan}, 
and not any compact totally geodesic spacelike surface in $\S_1^3$
is a spacelike slice. 
\end{rema}

\begin{rema} \label{refundido} The compactness of $S$ in Theorem \ref{tc1} cannot be relaxed to
completeness in general. In fact, a complete maximal surface in Lorentz-Minkowski spacetime $\L^ 3$ is a 
spacelike plane (see for instance \cite{Kobayashi}), and, clearly, not any spacelike plane is a 
spacelike slice in $\L^3$. However, complete spacelike surfaces in certain GRW spacetimes, which are
far from Lorentzian products, have been studied in \cite{Ro-Ru3}. In fact, in that reference the warping 
function is assumed non locally constant, i.e., there is no open rectangle $J\times F$ in $M$
such that, the restriction of the Lorentzian metric (\ref{Lorentzian_metric}) to $J\times F$ is a product.  
Such a GRW spacetime is called \emph{proper}. 
\end{rema}

\begin{rema}{\rm Taking into account \cite[Th. 3.1]{A-M},
\cite[Cor. 5.10]{Ca-Ro-Ru2}, the mean curvature $H_0$ of a compact CMC spacelike 
surface in a GRW spacetime, whose warping function satisfies $(\log f)'' \leq 0$, 
satisfies $\mathrm{(\;\widetilde I\;)}$ (really, the equality holds). Altenatively,
this fact can be deduced from the argument which follows. In fact, denote as previously 
by $\mathcal{F}$ a primitive function of $f$ and by $\mathcal{F}(t)$ its restriction on 
$S$. From (\ref{Laplacian}) we get $$\Delta \mathcal{F}(t)=-2f'(t)-2H_0f(t)\langle N,\partial_t\rangle.$$
Let $p_0$ and $p^0$ be the points of $S$ where $\mathcal{F}(t)$ attains its global minimum and
maximum values, respectively. We have $\nabla t(p_0)=\nabla t(p^0)=0$ and $\langle
N,\partial_t\rangle(p_0)=\langle N,\partial_t\rangle(p^0)=1$. On
the other hand, $\Delta\mathcal{F}(t)(p_0))\geq 0$ and
$\Delta \mathcal{F}(t)(p^0))\leq 0$ hold, and therefore the constant $H_0$ satisfies

$$\frac{-f'(t(p^0))}{\,f(t(p^0))}\,\leq \,H_0\,\leq\,\frac{-f'(t(p_0))}{\,f(t(p_0))}.$$

Now, it is enough to observe that the function $\mathcal{F}$
is increasing from its definition, and the function
$-f'/f$ is also increasing from our assumption, giving the
equality $H_0=-f'(t)/f(t)$. Therefore, for the case
of $3$-dimensional GRW spacetimes, the previously quoted result in \cite{A-M} follows
from Theorem \ref{tc1}}.

\vspace{1mm}

The same argument does not work for the case of noncompact complete spacelike surfaces. 
In fact, as an application of the classical generalized maximum principle, due to Omori \cite{Omori} 
and Yau \cite{Yau}, it was shown, \cite[Prop. 5.3]{Ro-Ru1}, that a complete spacelike 
surface of constant mean curvature in certain proper RW spacetimes with fiber $\mathbb{R}^2$,
warping function satisfying $(\log f)''\leq 0$ and which is contained between two spacelike 
slices must satisfy $\mathrm{(\;\widetilde I\;)}$.
\end{rema}

\begin{coro}\label{otro} Let $M$ be a 
GRW spacetime, whose warping function satisfies $(\log f)''\leq 0$
 and the Gauss curvature of its fiber
obeys $K^ F(\pi_F)\geq c \,f(\pi_I)^ 2$ for some constant $c>0$. If
$M$ admits a complete spacelike surface $S$ such that the
inequality $\mathrm{(\;\widetilde I\;)}$ holds on all $S$, then $M$ is spatially closed 
and $S$ is a spacelike slice.
\end{coro}

\noindent\emph{Proof.}  From (\ref{K2}), the Gauss curvature
of $S$ satisfies $K\geq c >0$. Therefore, $S$ is compact  from the classical Bonnet-Myers
theorem. Consequently, \cite[Prop. 3.2]{A-R-S1}, $F$ must be also compact. Finally, $S$ is a 
spacelike slice as a consequence of Theorem \ref{tc1}  according 
the assumption made on the warping function.
\hfill{$\Box$}

\begin{coro}\label{cc1} Let $M$ be a 
GRW spacetime whose warping function satisfies $(\log f)''\leq 0$. Let $S$ a compact spacelike surface in $M$
and put $\alpha=\min_S\frac{f'(t)^2}{f(t)^2}$. If the mean
curvature function $H$ of $S$ satisfies
$$H^2\leq \alpha,$$

\noindent then $S$ must be a spacelike slice. Moreover, there is no
compact spacelike surface in $M$ whose mean curvature function
satisfies $H^2<\alpha$.
\end{coro}

\begin{teor}\label{BM} Let $M$ be a 
GRW spacetime, whose warping function satisfies $(\log f)''\leq 0$
 and whose fiber $F$ has Gauss curvature
$K^F\geq 0$. Then, there is no complete spacelike surface $S$ in $M$ such
that for some positive constant $\epsilon$, the inequality
$$\frac{f'(t)^2}{f(t)^2}-H^2>\epsilon>0,$$
\noindent holds on all $S$.
\end{teor}
\emph{Proof.} Otherwise suppose there exists such a surface $S$. From (\ref{K2}), 
the Gauss curvature of $S$ satisfies $K\geq\epsilon >0$. Therefore, the Bonnet-Myers
theorem asserts that $S$ is compact. Consequently,
\cite[Prop. 3.2]{A-R-S1}, $F$ must be also compact. Moreover, as a
consequence of Theorem \ref{tc1},
$S$ is a spacelike slice with $H^2=\frac{f'(t)^2}{f(t)^2}$, which is a contradiction.
\hfill{$\Box$}

\begin{rema} As shown in \cite[Th. 6.20]{C-Ro-Ru-2}, the only complete 
spacelike surfaces of constant mean curvature in the 3-dimensional steady 
state spacetime that are bounded away from the infinity future and whose mean
curvature satisfies $H\leq 1$ are the spacelike slices. As application of the previous 
result, we have that there is no complete spacelike surface in the 3-dimensional
steady state spacetime, whose mean curvature function satisfies
$H^2<1-\epsilon$, for a constant $\epsilon$ such that $0<\epsilon<1$.
\end{rema}

\begin{rema}\label{npar} {\bf a)} The assumption on $K^F$ in Theorem
\ref{BM} is crucial to get the compactness of $S$ thanks to 
the classical Bonnet-Myers theorem. However, if the compactness 
of $S$ is directly assumed then, with no assumption 
on $K^F$, the following nonexistence result holds:  if
$M$ is a proper (resp. non-necessarily proper) GRW spacetime,
whose warping function satisfies $(\log f)''\leq 0$ (resp. $(\log
f)''<0$), then there is no compact spacelike surface $S$ in
$M$ such that its mean curvature function satisfies
$\frac{f'(t)^2}{f(t)^2}-H^2>0$ everywhere on $S$. In fact, from
(\ref{laplacian3}) we have that there exists $p\in S$ such
that $H^2(p)\geq \frac{f'(t(p))^2}{f(t(p))^2}$ (alternatively, $p$
may be taken as point where $\log f(t)$ attains its minimum
value). {\bf b)} Previous argument implies that when
inequality $\mathrm{(\;\widetilde I\;)}$  is assumed on a
compact spacelike surface $S$ in a GRW spacetime $M$ 
whose warping function satisfies $(\log f)'' \leq 0$, there is always 
some point of $S$ where the inequality for $H$ is indeed an 
equality. In the case $S$ is complete and noncompact  this is
not true, in general.
\end{rema}

We end this section with the case that the GRW spacetime $M$ has
constant sectional curvature $\bar c$. Note that this holds if and
only if its fiber $F$ has constant Gauss curvature $K^F=c$
(i.e., $M$ is a RW spacetime) and its 
warping function $f$ satisfies,

\begin{equation}\label{constant}
\frac{f''}{f}=\frac{c+(f')^2}{f^2}=\bar c,
\end{equation}
(see for instance \cite[Cor. 9.107]{Besse}). All the positive
solutions of (\ref{constant}) were found in \cite{A-R-S2}, as a
particular case of a much general result. Note that we have,

$$(\log f)''=\frac{c}{f^2}.$$
directly from (\ref{constant}). As shown in \cite[Table, cases 5,6]{A-R-S2},
if $\bar c \leq 0$ then $c<0$. Therefore, as a consequence of 
the previous results, we obtain the following extension of 
\cite[Cor. 5]{A-R-S2} for the case of spacelike surfaces.

\begin{teor}\label{RW}
The only compact spacelike surfaces $S$ in a RW spacetime 
with non-positive constant sectional curvature and whose mean
curvature function $H$ satisfies the inequality
$\mathrm{(\;\widetilde I\;)}$ on all $S$,  are the spacelike slices.
\end{teor}

\section{Uniqueness results of Calabi-Bernstein type}
Let $(F,g)$ be a $2$-dimensional compact Riemannian manifold and
let  $f : I \longrightarrow \R$ be a positive smooth function. For
each $u \in C^{\infty}(F)$ such that $u(F)\subset I$ we can
consider its graph $\Sigma_u=\{(u(p),p) \, : \, p\in F\}$ in the
$3$-dimensional GRW spacetime $M$ with base $(I,-dt^2)$, fiber
$(F,g)$ and warping function $f$. The graph of $u$ inherits a
metric, given by (\ref{induced_metric}), which is
Riemannian if and only if $u$ satisfies $\vert Du\vert<f(u)$
everywhere on $F$. Note that $t(u(p),p)=u(p)$ for any $p \in F$, and so, 
the functions $t$ on $\Sigma_u$ and  $u$ can be naturally identified. 
As an application of the previous results, we havee the following uniqueness 
results,

\begin{teor}\label{grafo1} Let $F$ be a $2$-dimensional compact Riemannian 
manifold and let $f:I \longrightarrow ]0,\infty[$ be a  smooth function such
that $(\log f)''\leq 0$. Then, the only
entire solutions $u: F\longrightarrow I$ of the inequality
$\mathrm{(I)}$ are the constant functions. Moreover, there is no entire
solution of the inequality assuming that $\mathrm{(I.1)}$ holds
strictly.
\end{teor}

\begin{teor}\label{grafo3} Let $F$ be a $2$-dimensional compact Riemannian manifold 
and let $f:I \longrightarrow ]0,\infty[$ be a  smooth function
with $(\log f)''\leq 0$  and such that
$\inf\, \frac{(f')^2}{f^2}=\alpha>0$. Then, the only entire
solutions $u: F\longrightarrow I$ of the inequality
$$H(u)^2\leq\alpha ,$$
$$\vert Du\vert<f(u),$$ are the constant functions. Moreover, there
is no entire solution if the first inequality holds strictly.
\end{teor}

\end{document}